\documentclass[a4paper,11pt,DIV=12,
abstract=on
]{scrartcl}
\usepackage{mathtools}
\mathtoolsset{showonlyrefs}
\usepackage{amsmath, amsthm, amssymb}
\usepackage{fontenc}
\usepackage[utf8]{inputenc}
\usepackage[english]{babel}
\usepackage{graphicx}
\usepackage{subcaption,xargs}
\usepackage{booktabs}
\usepackage{enumitem,listings}
\usepackage{environ,ifthen,xcolor}
\usepackage{hyperref}
\usepackage[
  defernumbers=true,
  backend=biber,
  style=numeric-comp,
  firstinits=true,
  natbib=true,
  maxcitenames=5,
  sortlocale=en_EN,
  url=true, 
  doi=true,
  eprint=true,
  maxnames=99
]{biblatex}
\setkomafont{sectioning}{\rmfamily\bfseries}
\setkomafont{title}{\rmfamily}

\newcommandx{\PT}{\operatorname{PT}}
\newcommandx{\M}{\mathcal{M}}
\newcommand{\tT}{\mathrm{T}}
\newcommand{\T}{\mathrm{T}}

\makeatother

\DeclareMathOperator*{\argmax}{arg\,max}

\title{Nonlocal Inpainting of Manifold-valued Data\\ on Finite Weighted Graphs}
\author{Ronny Bergmann%
	\thanks{Fachbereich für Mathematik,%
	Technische Universität Kaiserslautern, 67663 Kaiserslautern, Germany.\newline
	\texttt{bergmann@mathematik.uni-kl.de}}%
	\and%
	Daniel Tenbrinck%
	\thanks{Institute for Computational and Applied Mathematics,
	Westfälische Wilhelms-Universität Münster, 48149 Münster, Germany.\newline
	\texttt{daniel.tenbrinck@uni-muenster.de}}%
}
\hypersetup{pdfauthor={Ronny Bergmann, Daniel Tenbrinck},
	pdftitle={A Graph Framework for Manifold-Valued Data},
	pdfsubject={arXiv Preprint},%
	pdfcreator = {pdflatex and TextMate},%
	pdfkeywords={},
	plainpages=false, pdfstartview=FitH, pdfview=FitH, pdfpagemode=UseOutlines,%
	bookmarksnumbered=true,bookmarksopen=false,bookmarksopenlevel=0,%
	colorlinks=true,
		linkcolor=green!50!black,
		citecolor=red!50!black,
		urlcolor=blue!50!black%
}
\date{July 12, 2017}
\begin{document}
  \maketitle
  \begin{abstract}
  Recently, there has been a strong ambition to translate models and algorithms
  from traditional image processing to non-Euclidean domains, e.g., to
  manifold-valued data. While the task of denoising has been extensively
  studied in the last years, there was rarely an attempt to perform image
  inpainting on manifold-valued data. In this paper we present a nonlocal
  inpainting method for manifold-valued data given on a finite weighted graph.
  We introduce a new graph infinity-Laplace operator based on the idea of
  discrete minimizing Lipschitz extensions, which we use to formulate the
  inpainting problem as PDE on the graph. Furthermore, we derive an explicit
  numerical solving scheme, which we evaluate on two classes of synthetic
  manifold-valued images.
  \end{abstract}
  %
  %
  \section{Introduction}
  Variational methods and partial differential equations (PDEs) play a key role
  for both modeling and solving image processing tasks. When processing real
  world data certain information might be missing due to structural artifacts,
  occlusions, or damaged measurement devices. Reconstruction of missing image
  information is known as inpainting task and there exist various variational
  models to perform inpainting. One successful method from the literature is
  based on a discretization of the \(\infty\)-Laplace operator~\cite{CMS98}.
  This idea has been adapted to finite weighted graphs in~\cite{EDL12}.
  The graph model enables to perform local and nonlocal inpainting within the
  same framework based on the chosen graph construction. Nonlocal inpainting
  has the advantage of preserving structural features by using all available
  information in the given image instead of only local neighborhood values.
  
  With the technological progress in modern data sensors there is an emerging 
  field of processing non-Euclidean data. We concentrate our
  discussion in the following on manifold-valued data, i.e., each data value
  lies on a Riemannian manifold. Real examples for manifold-valued images are
  interferometric synthetic aperture radar (InSAR) imaging~\cite{MF98}, where
  the measured phase-valued data may be noisy and/or incomplete. Sphere-valued
  data appears, e.g., in directional analysis~\cite{MJ00}.
  Another application is diffusion tensor imaging
  (DT-MRI)~\cite{BML94}, where the diffusion tensors can be represented as
  \(3\times 3\) symmetric positive definite matrices, which also constitute a
  manifold. For such data, there were several variational methods and
  algorithms proposed to perform image processing tasks, see e.g.,
  \cite{BBSW16,BPS16,CS13,LSKC13,WDS14}. Recently, the authors
  generalized the graph \(p\)-Laplacian for manifold-valued data, \(1\leq p <
  \infty\), in~\cite{BT17} and derived an explicit as well as an semi-implicit
  iteration scheme for computing solutions to related partial difference equations on graphs as mimetic approximation of continuous PDEs. While the previous 
  work concentrated on denoising, the present work deals with the task of image
  inpainting of manifold-valued data. For this, we extend the already defined
  family of manifold-valued graph \(p\)-Laplacians by a new operator, namely
  the graph \(\infty\)-Laplacian for manifold valued data. We derive an
  explicit numerical scheme to solve the corresponding PDE and illustrate its
  capabilities by performing nonlocal inpainting of synthetic manifold-valued
  data.
  
  The remainder of this paper is organized as follows: In Sec.~\ref{sec:pre} we
  introduce the necessary notations of Riemannian manifolds, finite weighted
  graphs, and manifold-valued vertex functions. In Sec.~\ref{sec:infLapl} we
  introduce a new graph~\(\infty\)-Laplace operator for manifold-valued data
  based on the idea of discrete minimizing Lipschitz extensions. Furthermore,
  we derive an explicit numerical scheme to solve the corresponding
  PDE~\(\Delta_\infty f = 0\) with suitable boundary conditions. In
  Sec.~\ref{sec:num} we apply the proposed method to inpainting of synthetic 
  manifold-valued images. Finally, Sec.~\ref{sec:concl} concludes the paper.
  %
  %
  \section{Preliminaries}
  \label{sec:pre}
  In this section we first introduce the needed theory and notations on
  Riemannian manifolds in Sec.~\ref{subsec:manifolds} and introduce finite
  weighted graphs in Sec.~\ref{subsec:graphs}. We then combine both concepts to
  introduce vertex functions and tangential edge functions needed for the
  remainder of this paper in Sec.~\ref{subsec:functions}. For further details we
  refer to \cite{BT17}.

  %
    \subsection{Riemannian manifolds}
  \label{subsec:manifolds}
  For a detailed introduction to functions on Riemannian manifolds we refer to,
  e.g., \cite{AMS08,Jost11}. The values of the given data lie in a
  \emph{complete, connected, \(m\)-di\-men\-sion\-al Riemannian manifold}~\(\mathcal
  M\) with Riemannian metric~\(\langle \cdot,\cdot\rangle_x \colon \T_x\mathcal
  M\times \T_x\mathcal M\to\mathbb R\), where~\(\T_x\mathcal M\) is the
  tangent space at \(x \in \M\). In every tangent space \(\T_x\mathcal
  M\) the metric induces a norm, which we denote by~\(\lVert\cdot\rVert_x\).
  The disjoint union of all tangent spaces is called the \emph{tangent
  bundle}~\(\T\mathcal M \coloneqq \dot\cup_{x\in\mathcal M}\T_x\mathcal M\).
  Two points \(x,y\in\mathcal M\) can be joined by a (not necessarily unique)
  shortest curve~\(\gamma_{\overset{\frown}{x,y}}\colon [0,L] \to \mathcal M\),
  where \(L\) is its length. This generalizes the idea of shortest paths from
  the Euclidean space~\(\mathcal M=\mathbb R^m\), i.e., straight lines, to a
  manifold and induces the geodesic distance denoted~\(d_{\mathcal M}\colon
  \mathcal M\times \mathcal M\to\mathbb R^+\). A curve~\(\gamma\) can be
  reparametrized such that derivative vector field~\(\dot\gamma(t) \coloneqq
  \frac{d}{dt}\gamma(t)\in T_{\gamma(t)}\mathcal M\) has constant norm,
  i.e.,~\(\lVert \dot\gamma_{\overset{\frown}{x,y}}(t)
  \rVert_{\gamma_{\overset{\frown}{x,y}}(t)} = 1\), \(t\in[0,L]\). The
  corresponding curve then has unit speed. We employ another notation of a
  geodesic, namely \(\gamma_{x,\xi}\), \(\xi\in\T_x\mathcal M\), \(x\in\mathcal
  M\), which denotes the locally unique geodesic fulfilling \(\gamma_{x,\xi}(0)
  = x\) and~\(\dot\gamma_{x,\xi}(0) = \xi\in \T_x\mathcal M\). This is unique
  due to the Hopf--Rinow Theorem, cf.~\cite[Theorem 1.7.1]{Jost11}. We further
  introduce the \emph{exponential map}~\(\exp_x\colon \T_x\mathcal M\to\mathcal
  M\) as~\(\exp_x(\xi) = \gamma_{x,\xi}(1)\). Let~\( r_x\in\mathbb R^+\) denote
  the injectivity radius, i.e., the largest radius such that \(\exp_x\) is
  injective for all \(\xi\) with~\(\lVert\xi\rVert_x < r_x\). Furthermore, let
  \[ \mathcal D_x \coloneqq \bigl\{ y\in\mathcal M : y = \exp_x\xi, \text{ for
  some }\xi\in T_x\mathcal M \text{ with } \lVert\xi\rVert_x<r_x\bigr\}\ . \]
  Then the inverse map~\(\log_x\colon \mathcal D_x\to \T_x\mathcal M\) is
  called the~\emph{logarithmic map} and maps a
  point~\(y=\gamma_{x,\xi}(1)\in\mathcal D_x\) to~\(\xi\).

  %
    \subsection{Finite weighted graphs}%
  \label{subsec:graphs}
  Finite weighted graphs allow to model relations between arbitrary discrete
  data. Both local and nonlocal methods can be unified within the graph
  framework by using different graph construction methods: spatial vicinity for
  local methods like finite differences, and feature similarity for nonlocal
  methods. A finite weighted graph \(G=(V,E,w)\) consists of a finite
  set of indices~\(V=\{1,\ldots,n\}\), \(n\in \mathbb N\), denoting the
  vertices, a set of directed edges \(E\subset V\times V\) connecting a subset
  of vertices, and a nonnegative weight function~$w \colon E \rightarrow
  \mathbb{R}^+$ defined on the edges of a graph. For an edge \((u,v)\in E\),
  \(u,v\in V\) the node \(u\) is the start node, while \(v\) is the end node.
  We also denote this relationship by \(v\sim u\). Furthermore, the weight
  function \(w\) can be extended to all \(V\times V\) by setting \(w(u,v) = 0\)
  when \(v\not\sim u\). The neighborhood \(\mathcal N(u) \coloneqq \{ v\in V :
  v\sim u \}\) is the set of adjacent nodes.

  %
    \subsection{Manifold-valued vertex functions and tangential edge functions.}
    \label{subsec:functions}
    The functions of main interest in this work are manifold-valued vertex
    functions, which are defined as
    \[
      f\colon V\to\mathcal M,\quad
      u\mapsto f(u),
    \]
    The range of the vertex
    function~\(f\) is the Riemannian manifold \(\mathcal M\).
    We denote the set of admissible vertex functions by \(\mathcal H(V;\mathcal M) \coloneqq \{ f\colon V\to\mathcal M \} \).
    This set can be equipped with a metric given by
    \begin{equation*}
       d_{\mathcal H(V;\mathcal M)}(f,g)
       \ \coloneqq \
        \biggl(
          \sum_{u\in V} d_{\mathcal M}^2(f(u),g(u))
          \biggr)^{\frac{1}{2}},
          \qquad f,g\in\mathcal H(V; \mathcal M).
    \end{equation*}
    Furthermore, we need the notion of a tangential vertex function. The space
    \(\mathcal H(V;\T\mathcal M)\) consists of all functions \(H\colon
    V\to\T\mathcal M\), i.e., for each \(u\in V\) there exists a value
    \(H(u)\in\T_x\mathcal M\) for some \(x \in \mathcal M\).
  %
  \section{Methods}\label{sec:infLapl}
  In this section we generalize the \(\infty\)-Laplacian from the real-valued,
  continuous case to the manifold-valued setting on graphs. We
  discuss discretizations of the \(\infty\)-Laplacian both for
  real-valued functions on bounded open sets and on graphs in
  Sec.~\ref{subsec:discretizations}. We generalize these to manifold-valued
  functions on graphs in Sec.~\ref{subsec:manInfLaplace} and state a
  corresponding numerical scheme in Sec.~\ref{subsec:numScheme}.
  \subsection{Discretizations of the \(\infty\)-Laplace
  operator}\label{subsec:discretizations} Let \(\Omega \subset \mathbb{R}^d\)
  be a bounded, open set and let \(f \colon \Omega \rightarrow \mathbb{R}\) be
  a smooth function. Following~\cite{CEG01} the infinity Laplacian of \(f\) at
  \(x \in \Omega\) can be defined as
  \begin{equation}
  \label{eq:inftyLaplacian}
    \Delta_\infty f(x) \ = \
    \bigl(
    (\nabla f)^\tT\Delta f\nabla f
    \bigr)(x)
    \ = \  
      \sum\limits_{j=1}^d\sum\limits_{k=1}^d
        \frac{\partial f}{\partial x_j}
        \frac{\partial f}{\partial x_k}
        \frac{\partial^2 f}{\partial x_jx_k} (x).
  \end{equation}
  As discussed above, this operator is not only interesting in theory, but also
  has applications in image processing~\cite{CMS98}, e.g., for image
  interpolation and inpainting. Oberman discussed in~\cite{Obe04} different
  possibilities for a consistent discretization scheme of the infinity
  Laplacian defined in \eqref{eq:inftyLaplacian}. One basic observation is that 
  the operator can be well approximated by the maximum and minimum values of the
  function in a local \(\varepsilon\)-ball neighborhood, i.e.,
  \begin{equation}
  \label{eq:discretizationBall}
    \Delta_\infty f(x) \ = \ 
    \frac{1}{\varepsilon^2} \left( \min_{y \in B_\varepsilon(x)} f(y)
    + \max_{y \in B_\varepsilon(x)} f(y) - 2f(x) \right) \, 
    + \, \mathcal{O}(\varepsilon^2).
  \end{equation}
  The approximation in~\eqref{eq:discretizationBall} has inspired Elmoataz et
  al.~\cite{EDL12} to propose a definition of a discrete graph
  \(\infty\)-Laplacian operator for real-valued vertex functions, i.e.,
  \begin{equation}
    \label{eq:graphInf}
    \begin{split}
	  \Delta_\infty f(u) \ &=\ 
    \max_{v\sim u} \, \lvert \max(\sqrt{w(u,v)} (f(v) - f(u)),0) \rvert
    \\ \: &\hspace{1cm}- \ 
    \max_{v\sim u} \, \lvert \min(\sqrt{w(u,v)} (f(v) - f(u)),0) \rvert
    \end{split}
    \end{equation}
  Furthermore, Oberman uses in \cite{Obe04} the well-known relationship between
  solutions of the homogeneous infinity Laplace equation \(-\Delta_\infty f =
  0\) and absolutely minimizing Lipschitz extensions to derive a numerical
  scheme based on the idea of minimizing the discrete Lipschitz constant in a
  neighborhood.
  \subsection{The graph \(\infty\)-Laplacian for manifold-valued data}\label{subsec:manInfLaplace}
    Instead of following the approach proposed by Elmoataz et al.~in
    \cite{EDL12} we propose a new graph \(\infty\)-Laplace operator for
    manifold valued functions based on the idea of computing discrete minimal
    Lipschitz extensions, i.e., for a vertex function \(f \in \mathcal H(V;
    \mathcal M)\) we define the graph \(\infty\)-Laplacian for manifold valued
    data \( \Delta_\infty \colon \mathcal H(V; \mathcal M) \rightarrow H(V; T
    \mathcal M\)) in a vertex \(u \in V\) as
    \begin{equation}
      \label{eq:InfLaplacianNew}
      \Delta_\infty f(u) \ \coloneqq \ 
        \frac{\left(
          \sqrt{w(u,v_1^*)}\log_{f(u)} f(v_1^*)
          \: + \: \sqrt{w(u,v_2^*)} \log_{f(u)} f(v_2^*)
        \right)}{\sqrt{w(u,v_1^*)} \, + \sqrt{w(u,v_2^*)}} 
  \end{equation}
  for which the designated neighbors \(v_1^*, v_2^* \in \mathcal N(u)\) are
  characterized by maximizing the discrete Lipschitz constant in the local
  tangential plane \(T_{f(u)}\mathcal M\) among all neighbors, i.e.,
  \begin{equation*}
    (v_1^*, v_2^*) \ = \!\!\!
  \argmax_{(v_1, v_2) \in \mathcal N^2(u)} 
  \Bigl\lVert
    \sqrt{w(u,v_1)}\log_{f(u)} f(v_1) - \sqrt{w(u,v_2)}\log_{f(u)} f(v_2)
  \Bigr\rVert_{f(u)}
  \end{equation*}
  By means of the proposed operator in \eqref{eq:InfLaplacianNew} we are
  interested in solving discrete interpolation problems on graphs for
  manifold-valued data. Let \(U \subset V\) be a subset of vertices of the
  finite weighted graph \(G = (V,E,w)\) and let \(f \colon V/U \rightarrow
  \mathcal M\) be a given vertex function on the complement of \(U\). The
  interpolation task now consists in computing values of \(f\) on \(\mathcal
  M\) for vertices \(u \in U\) in which \(f\) is unknown. For this we solve 
  the following PDE on a graph based on the proposed operator
  in~\eqref{eq:InfLaplacianNew} with given boundary conditions:
    \begin{equation}
    \label{eq:PDE}
    \begin{cases}
	    \ \Delta_\infty f(u) \, = \ 0 \qquad &\text{ for all } u \in U, \\
	    \ f(u) \, = \ g(u) &\text{ for all } u\in V/U.
    \end{cases}
  \end{equation}
  \subsection{Numerical iteration scheme}
    \label{subsec:numScheme}
    In order to numerically solve the PDE in~\eqref{eq:PDE} on a finite
    weighted graph we introduce an artificial time variable \(t\) and derive a
    related parabolic PDE, i.e.,
    \begin{equation}
      \label{eq:parabolicEq}
      \begin{cases}
        \ \frac{\partial f}{\partial t}(u,t)
        \, = \ \Delta_\infty f(u,t) \qquad &\text{ for all } u \in U,
          \, t \in (0, \infty),\\
          \ f(u,0) \, = \ f_0(u) &\text{ for all } u\in U,\\
        \ f(u,t) \, = \ g(u,t) &\text{ for all } u\in V/U, t \in [0, \infty).
      \end{cases}
    \end{equation}
    We propose an explicit Euler time discretization scheme with sufficiently small
    time step size \( \tau > 0\) to iteratively solve~\eqref{eq:parabolicEq}. 
    Note that we proposed a similar explicit scheme for the computation of
    solutions of the graph \(p\)-Laplacian operator for manifold-valued 
    data in~\cite{BT17}. Using the notation \(f_k(u) \coloneqq f(u,k\tau)\), i.e. we discretize the time \(t\) in steps \(k\tau\), \(k\in\mathbb N\), an 
    update for the vertex function \(f\) can be computed by
    \begin{equation}\label{eq:explInfLaplace}
      f_{k+1}(u) \ = \ \exp_{f_k(u)}
      \bigl(
      \tau
      \Delta_\infty f_k(u)
      \bigr),
    \end{equation}
    for which the graph \(\infty\)-Laplacian is defined
    in~\eqref{eq:InfLaplacianNew} above.
  %
  %
  \section{Numerical examples}\label{sec:num}
  We first describe our graph construction and details of our inpainting
  algorithm in Sec.~\ref{ssec:graphConstr}. We then consider two synthetic 
  examples of manifold-valued data in Sec.~\ref{ssec:results}, namely 
  directional image data and an image consisting of symmetric positive matrices 
  of size~\(2\times2\).
  
  \subsection{Graph construction and inpainting algorithm}\label{ssec:graphConstr}
    We construct a nonlocal graph from a given manifold-valued image
    \(f\in\mathcal M^{m,n}\) as follows: we consider patches
    \(q_{i,j}\in\mathcal M^{2p+1,2p+1}\) of size \(2p+1\) around each pixel \((i,j), i\in\{1,\ldots,m\}, j\in \{1,\ldots,n\}\)
    with periodic boundary conditions. We denote for two patches \(q_{i,j}\)
    and \(q_{i',j'}\) the set \(I
    \subset\{1,\ldots,2p+1\}\times\{1,\ldots,2p+1\}\) as pixels that are known
    in both patches and compute~\(d_{i,j,i',j'} \coloneqq \frac{1}{\lvert I
    \rvert} d_{\mathcal H(I ;\mathcal M)}(q_{i,j},q_{i',j'})\).
    For~\((i,j)\in\{1,\ldots,n\}\times\{1,\ldots,m\}\) we introduce edges to
    the pixels \((i',j')\) with the \(k\) smallest patch distances. We define
    the weight function as~\(w(u,v) = \mathrm{e}^{-d^2_{i,j,i',j'} /
    \sigma^2}\), where \(u=(i,j),v=(i',j')\in\mathcal G\) with \(v\sim u\). For
    computational efficiency we further introduce a search window size
    \(r\in\mathbb N\), i.e., similar patches are considered within a window of
    size \(2r+1\) around \((i,j)\) around \(i,j)\), a construction that is for example also used for non-local means~\cite{BCM05}.

    We then solve the iterative scheme~\eqref{eq:parabolicEq} with
    \(\tau=\frac{1}{10}\) in~\eqref{eq:explInfLaplace} on all pixels \((i,j)\)
    that where initially unknown. We start with all border pixels,
    i.e.,~unknown pixels with at least one known local neighbor pixel. We stop
    our scheme if the relative change between two iterations falls below
    \(\epsilon = 10^{-7}\) or after a maximum of \(t=1000\) iterations. We then
    add all now border pixel to the active set we solve the equation on and
    reinitialize our iterative scheme~\eqref{eq:explInfLaplace}. We iterate
    this algorithm until all unknown pixel have been border pixel and are hence
    now known. Our algorithm is implemented in MathWorks MATLAB employing the
    Manifold-valued Image Restoration Toolbox (MVIRT)\footnote{open source, %
    \href{http://www.mathematik.uni-kl.de/imagepro/members/bergmann/mvirt/}%
    {www.mathematik.uni-kl.de/imagepro/members/bergmann/mvirt/}}.
  \subsection{Inpainting of directional data}
  \label{ssec:results}
  \begin{figure}[t]\centering
    \begin{subfigure}[t]{.49\textwidth}
      \centering
      \includegraphics[width=.9\textwidth]{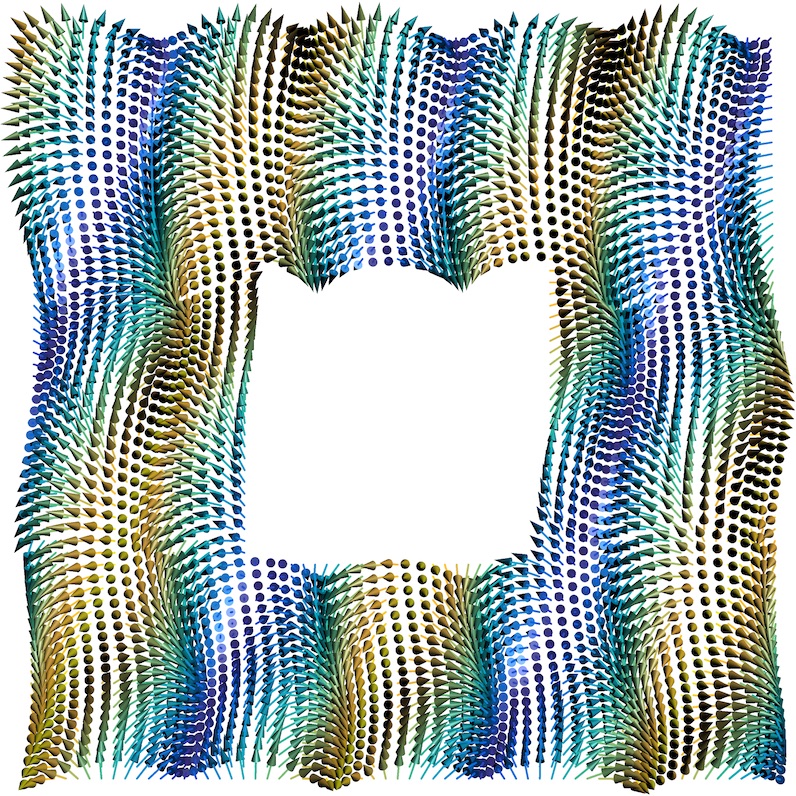}
      \vspace{-.33\baselineskip}
      \caption{Lossy Input data.}
      \label{subfig:S2:lossy}
    \end{subfigure}
    \begin{subfigure}[t]{.49\textwidth}
      \centering
      \includegraphics[width=.9\textwidth]{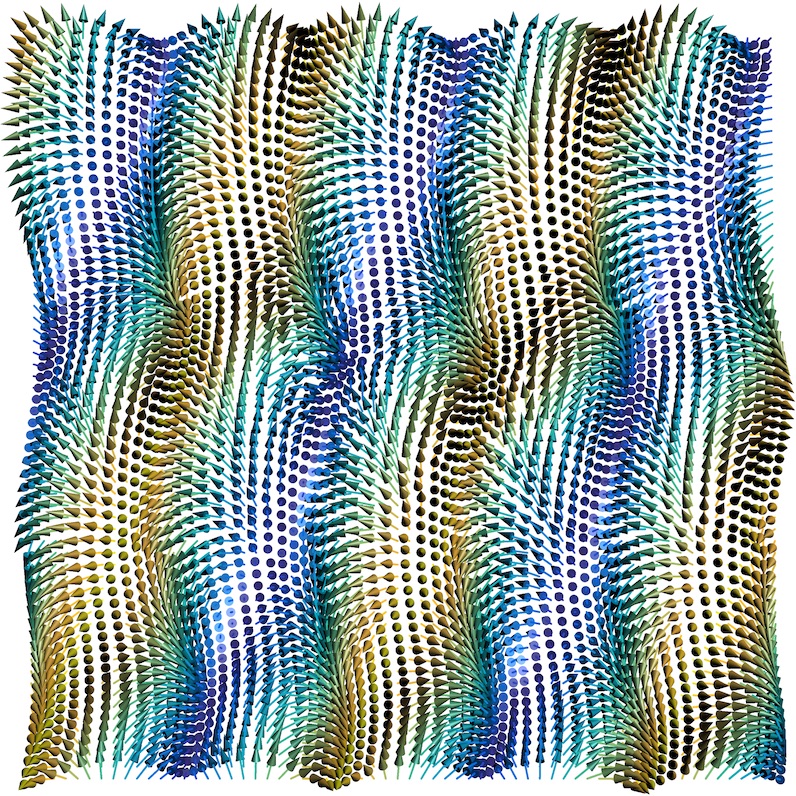}
      \vspace{-.33\baselineskip}
      \caption[]{Reconstruction,\\\(k=25\), \(p=12\).} \label{subfig:S2:k25p12}
    \end{subfigure}
    \vspace{-.5\baselineskip}
    \caption{Reconstruction of directional
      data~\(f\in(\mathbb S^2)^{64\times 64}\) from lossy given data, i.\,e.~
      in~(\subref{subfig:S2:lossy}) the original data is shown with missing data (the center).
      Using the \(k=25\) most similar patches and a patch radius of \(p=12\) leads to a reasonable reconstruction in~(\subref{subfig:S2:k25p12}).}
  \end{figure}
  We investigate the presented algorithm for artificial manifold-valued data: first, let \(\mathcal M = \mathbb S^2\) be
  the unit sphere, i.e., our data items are directions in \(\mathbb R^3\).
  Its data items are drawn as small three-dimensional arrows color-encoded by
  elevation, i.e.,~the south pole is blue (dark), the north pole yellow
  (bright). The periodicity of the data is slightly obstructed by two vertical
  and horizontal discontinuities of jump height~\(\frac{\pi}{16}\) dividing the
  image into nine parts. The input data is given by the lossy data,
  cf.~Fig.~\ref{subfig:S2:lossy}). We set the search window to~\(r=32\), i.e.,~global comparison of patches in the graph construction from Sec.~\ref{ssec:graphConstr}. Using~\(k=25\) 
  most similar patches and a patch radius of \(p=12\), the iterative scheme~\eqref{eq:explInfLaplace} yields Fig.~\ref{subfig:S2:k25p12}). The proposed methods finds a reasonable interpolation in the missing pixels.

  \subsection{Inpainting of symmetric positive definite matrices}
    \begin{figure}[t]\centering
      \begin{subfigure}[t]{.325\textwidth}
        \centering
        \includegraphics[width=.975\textwidth]{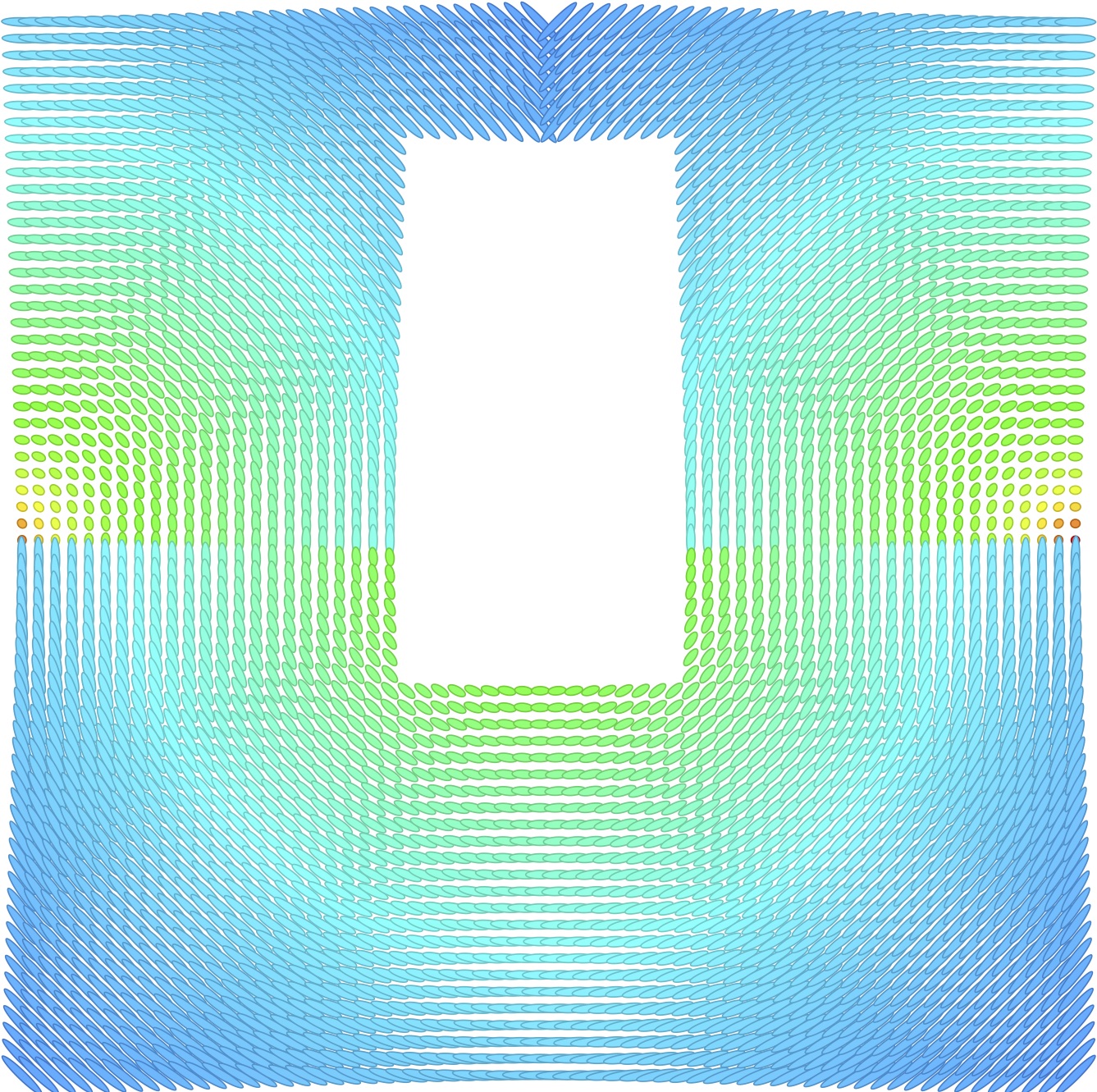}
      \vspace{-.33\baselineskip}
        \caption{Input data with missing rectangle.}
        \label{subfig:SPD2:lossy}
      \end{subfigure}
      \begin{subfigure}[t]{.325\textwidth}
        \centering%
        \includegraphics[width=.975\textwidth]{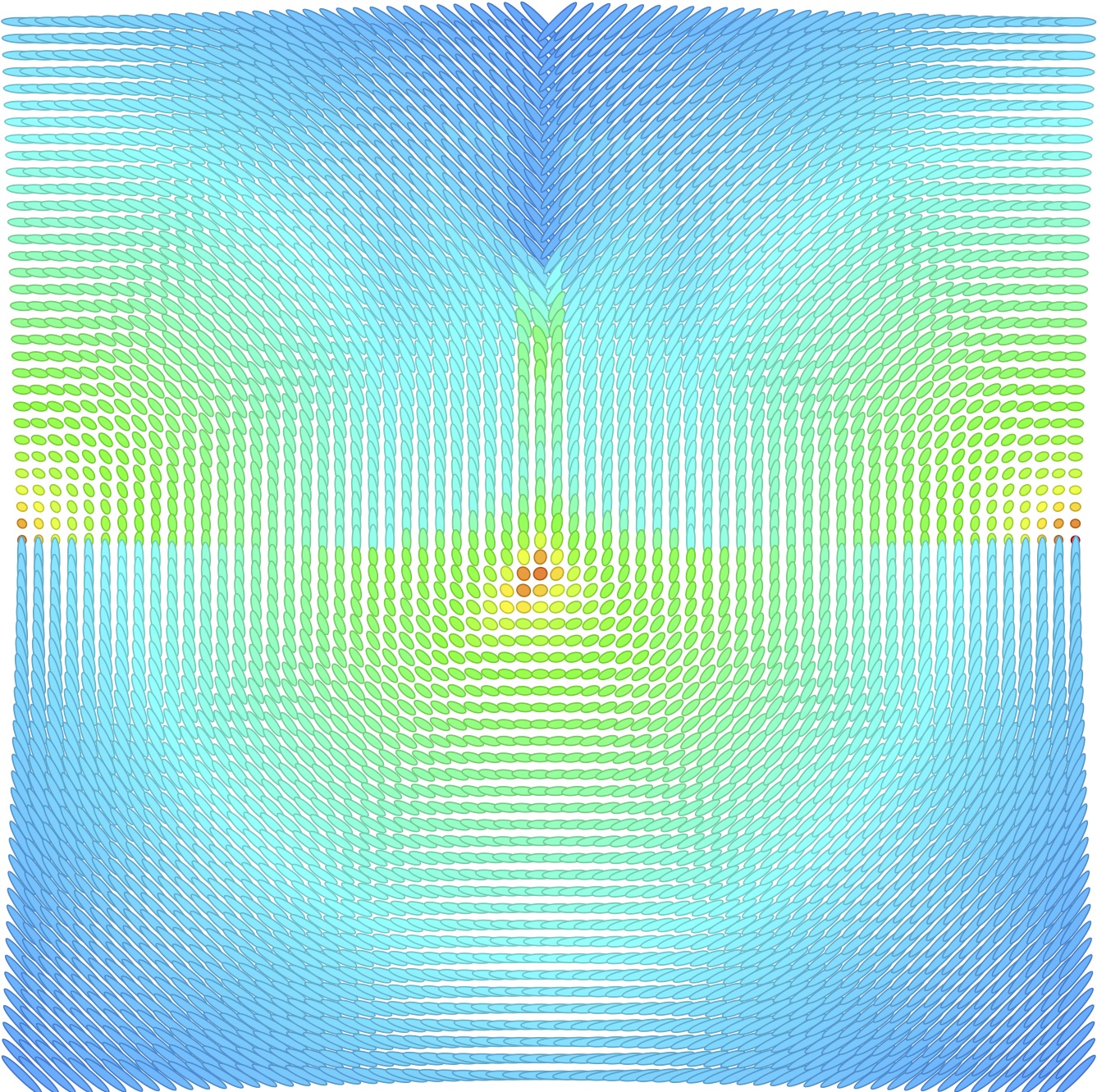}
      \vspace{-.33\baselineskip}
        \caption{Reconstruction, \(k=5\), \(p=6\).}
        \label{subfig:SPD2:recon1}
      \end{subfigure}
      \begin{subfigure}[t]{.325\textwidth}
        \centering%
        \includegraphics[width=.975\textwidth]{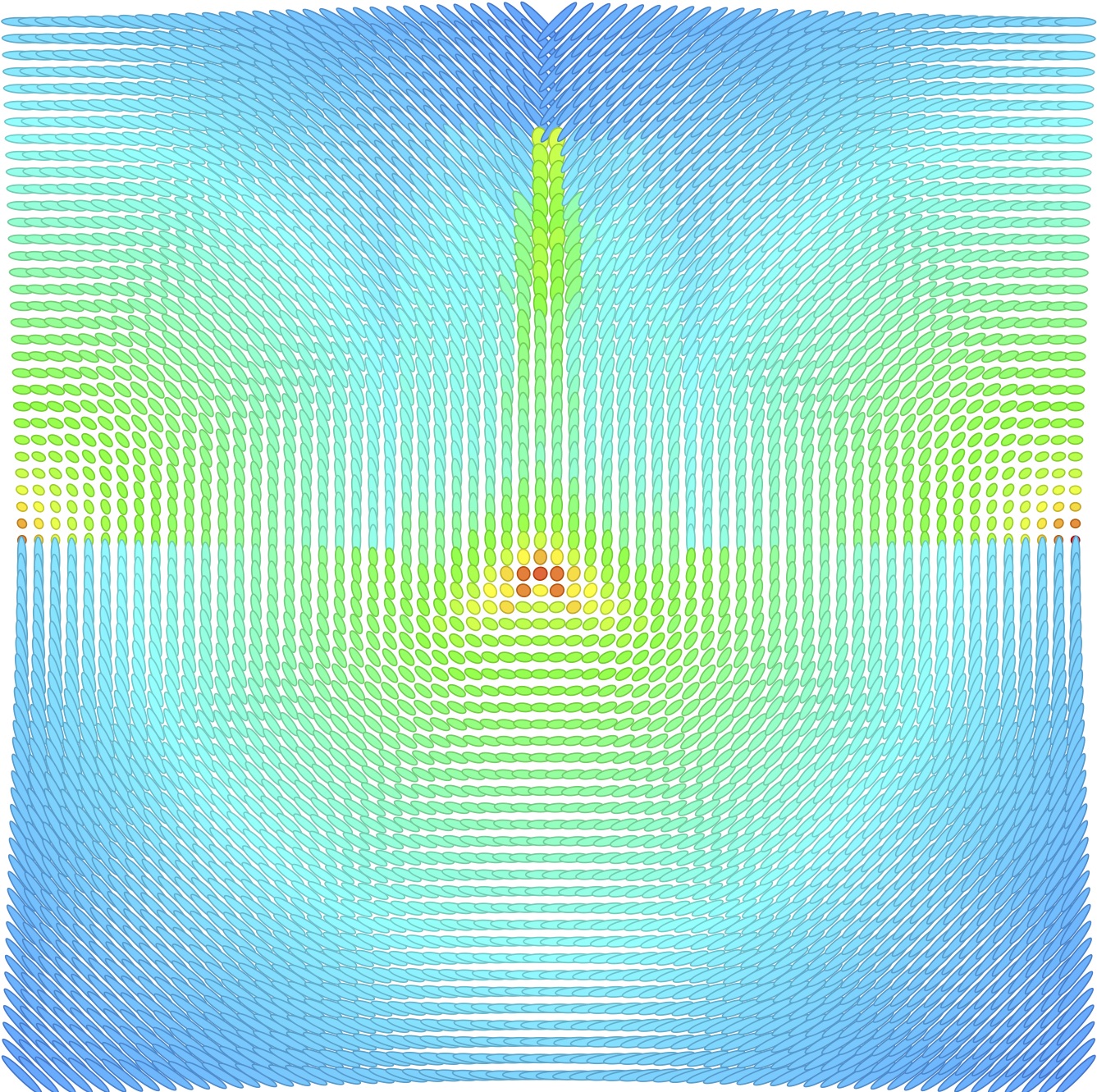}
      \vspace{-.33\baselineskip}
        \caption{Reconstruction, \(k=25\), \(p=6\).}
        \label{subfig:SPD2:recon2}
      \end{subfigure}
      \caption{
        Reconstruction of s.p.d.~matrices \(f\in(\mathcal
        P(2))^{64\times 64}\) of (\subref{subfig:SPD2:lossy}) a lossy image.
        Increasing the number of neighbors from \(k=5\)
        in~(\subref{subfig:SPD2:recon1}) to \(k=25\)
        in~(\subref{subfig:SPD2:recon2}) broadens the center feature and smoothens the discontinuity in the center.
      }
    \end{figure}
  As a second example we consider an image of symmetric positive definite (s.p.d.) matrices from \(\mathbb R^{2\times 2}\), i.e., \(\mathcal M = \mathcal P(2)\). These can be illustrated as ellipses using their eigenvectors as main axes and their eigenvalues as their lengths and the geodesic anisotropy index~\cite{MB06} in the hue colormap. The input data of \(64\times64\) pixel is missing a rectangular area, cf.~Fig.~\ref{subfig:SPD2:lossy}). We set again \(r=32\) for a global comparison. Choosing~\(k=5\), \(p=6\) yields a first inpainting result shown in Fig.~\ref{subfig:SPD2:recon1}) which preserves the discontinuity line in the center and introduces a red area within the center bottom circular structure. Increasing the nonlocal neighborhood to \(k=25\), cf.~in Fig.~\ref{subfig:SPD2:recon2}) broadens the red center feature and the discontinuity gets smoothed along the center vertical line.
  %
  %
  \section{Conclusion}\label{sec:concl}
  In this paper we introduced the graph \(\infty\)-Laplacian operator for manifold-valued functions by generalizing a reformulation of the \(\infty\)-Laplacian for real-valued functions and using discrete minimizing Lipschitz extensions. To the best of our knowledge, this generalization induced by our definition is even new for the vector-valued \(\infty\)-Laplacian on images. This case is included within this framework by setting e.g.~\(\mathcal M = \mathbb R^3\) for color images.
We further derived an explicit numerical scheme to solve the related parabolic PDE on a finite weighted graph.
First numerical examples using a nonlocal graph construction with patch-based similarity measures demonstrate the capabilities and performance of the inpainting algorithm applied to manifold-valued images.

Despite an analytic investigation of the convergence of the presented scheme, future work includes further development of numerical algorithms, as well as properties of the \(\infty\)-Laplacian for manifold-valued vertex functions on graphs.
\printbibliography
\end{document}